\documentclass[12pt]{article}
\usepackage[utf8]{inputenc}

\newtheorem{lem}{Lemma}
\newtheorem{thm}{Theorem}

\newtheorem{conj}{Conjecture}
\newcommand{\pf}{{ \noindent \bf Proof: }}
\newcommand{\qed}{\hfill $\Box$}
\usepackage{amssymb}
\begin{document}

\title{Subdivisions of maximal 3-degenerate graphs of order $d+1$ in graphs of minimum
degree $d$}

\author{Ajit A. Diwan\\
Department of Computer Science and Engineering,\\
Indian Institute of Technology Bombay,\\
Mumbai 400076, India.\\
\texttt{email}: aad@cse.iitb.ac.in
}

\maketitle

\abstract{ We prove that every graph of minimum degree at least $d \ge 1$ contains a
subdivision of \emph{some} maximal 3-degenerate graph of order $d+1$. This
generalizes the classic results of Dirac ($d=3$) and Pelikán ($d=4$). We conjecture
that for any planar maximal 3-degenerate graph $H$ of order $d+1$ and any graph
$G$ of minimum degree at least $d$, $G$ contains a subdivision of $H$. We verify this 
in the case $H$ is $P_6^3$ and $P_7^3$. }

\section{Introduction}

A classic result of Dirac~\cite{D} states that every graph of minimum degree at least
3 contains a subdivision of $K_4$.  Pelikán~\cite{P} proved that every graph of
minimum degree at least 4 contains a subdivision of $K_5^-$, the graph obtained
by deleting an edge from $K_5$.  In general, Mader~\cite{M1} first showed that there exists
a function $f(k)$ such that every graph of minimum degree at least $f(k)$ contains a subdivision of $K_k$. Bollobás and Thomason~\cite{BT} showed that $f(k)$ is 
$O(k^2)$ and this is best possible.

We consider the question in the other direction.  For a given integer
$d \ge 1$, for what graphs $H$ is it true that every graph of minimum degree at least $d$
contains a subdivision of $H$?  Clearly, any such graph $H$ can have at most $d+1$ 
vertices, since $K_{d+1}$ has minimum degree $d$. We consider graphs $H$ of order
exactly $d+1$.  We call a graph $H$ \emph{good} if every graph of minimum degree
at least $|H|-1$ contains a subdivision of $H$. For $1 \le d \le 3$, it follows that 
$K_{d+1}$ is good. Since there are planar graphs of minimum degree 4, $K_5$ 
is not good,  but Pelikán's theorem implies that $K_5^-$ is good.  We are interested in
finding the maximal good graphs. Mader~\cite{M2} showed that every graph of 
minimum degree at least $d \ge 2$ contains a  pair of adjacent vertices with $d$ internally 
disjoint paths between them.  This implies that the graph $K_2 \vee \overline{K_{d-1}}$, 
consisting of $d-1$ triangles that share a common edge, is good.  However, this graph 
has only $2d-1$ edges and is not a maximal good graph even for $d = 3$. Turner~\cite{T}
showed that the wheel $W_d = C_d \vee K_1$ is good, for all $d \ge 3$, but again
this has size $2d$ and is not a maximal good graph for $d = 4$.

Our main result is that every graph of minimum degree at least $d \ge 2$ contains a 
subdivision of \emph{some} graph $H$ of order $d+1$ and size $3d-3$. For $d = 3,4$ 
this implies the theorems of Dirac and Pelikán, respectively, since $K_4$ and $K_5^-$ 
are the only possible such graphs. Further, for $d = 5$, this is the maximum possible 
number of edges in a good graph, since there exist planar graphs of minimum degree 5.
We are unable to prove that any specific graph $H$ of order $d+1$ and size
$3d-3$ is good, for general $d$, but we can say something more about the structure of 
the graph $H$. 
We show that $H$ can be chosen to be 3-degenerate, that is, every subgraph of $H$ 
contains a vertex of degree at most 3.  We conjecture that every planar 3-degenerate graph of order $d+1$ and size $3d-3$ is good. We prove this for two specific graphs
$P_6^3$ and $P_7^3$. A weaker conjecture would be that $P_n^3$ is good for
all $n \ge 2$.

\section{Notation}

All graphs considered are undirected, finite and simple. The vertex set of a graph 
$G$ is denoted by $V(G)$ and the edge set by $E(G)$. The order of a graph $G$ is
$|V(G)|$ and $|E(G)|$ is its size. The subset of vertices adjacent to 
a vertex $v \in V(G)$ in a graph $G$ is denoted by $N_G(v)$ and $d_G(v) = |N_G(v)|$
is the degree of the vertex $v$. If $S \subset V(G) \cup E(G)$, $G-S$ is the graph obtained 
from $G$ by deleting all vertices and edges in $S$ and also edges incident with vertices
in $S$. A graph $H$ is a subgraph of a graph $G$ if $V(H) \subseteq V(G)$ and
$E(H) \subseteq E(G)$.  If $S \subset V(G)$, $G-(V(G)\setminus S)$ is the subgraph
of $G$ induced by $S$. The union of two subgraphs $H_1,H_2$ of a graph $G$ is
the subgraph with vertex set $V(H_1) \cup V(H_2)$ and edge set $E(H_1) \cup E(H_2)$.

A path $P$ in a graph $G$ is a sequence of distinct vertices $v_0,\ldots,v_{l}$ such
that $v_iv_{i+1}$ is an edge in $G$ for $0 \le i < l$. We say $P$ is a $v_0$--$v_l$ path
that joins $v_0$ to $v_l$. The vertices $v_0,v_l$ are the endpoints of $P$ and 
$\{v_1,\ldots,v_{l-1}\}$ are the internal vertices of $P$.  The set of internal vertices
of $P$ is denoted $I(P)$.  We will also consider $P$ to be a subgraph of $G$ with
vertex set $\{v_0,\ldots,v_l\}$ and edge set $\{v_iv_{i+1} \mid 0 \le i < l\}$.
A path $P$ is said to be an $A$--$B$ path in $G$, for
$A,B \subseteq V(G)$, if  $P$ joins a vertex  in $A$ to a vertex  in $B$ and
$I(P) \cap (A \cup B) = \emptyset$. An $A$--$B$ path is also said to join $A$ to $B$.
A set of paths $\mathcal{P}$ in $G$ is said to be internally disjoint if for any two distinct
paths $P,Q \in \mathcal{P}$, $I(P) \cap I(Q) = \emptyset$.  If $\mathcal{P}$
is a set of internally disjoint paths, let $I(\mathcal{P}) = \bigcup_{P \in \mathcal{P}}
\ I(P)$ be the set of internal vertices of $\mathcal{P}$. If $A \subset V(G)$ and
$u \in V(G) \setminus A$, a $u$--$A$ fan is a set of  internally disjoint $u$--$A$
paths having distinct endpoints in $A$. 

A graph $G$ is said to contain a subdivision of a graph $H$ if there exists a subset
$B(H) \subseteq V(G)$ of vertices and a set $\mathcal{P}$ of  internally
disjoint $B(H)$--$B(H)$ paths in $G$ such that:
\begin{enumerate}
\item
There exist bijections $f : V(H) \rightarrow B(H)$ and $g : E(H)
\rightarrow \mathcal{P}$.
\item
If $uv \in E(H)$ then $g(uv)$ is an $f(u)$--$f(v)$ path in $G$.
\end{enumerate}

We call the subgraph of $G$ formed by the union of the paths in $\mathcal{P}$
a subdivision of $H$ and denote it $\mathcal{T}(H)$. The vertex $f(v) \in 
V(\mathcal{T}(H))$ is said to correspond to the vertex $v \in V(H)$.

An ordered clique in a graph $G$ is a complete subgraph of $G$ together 
with a total ordering imposed on the vertices in the complete subgraph.

Let $G$ be a graph and $K$ an ordered clique in $G$. Let $u_1, u_2, \ldots, u_t$ be 
a sequence of vertices in $V(G) \setminus V(K)$ and $n_1, n_2, \ldots, n_t$ a 
sequence of positive integers. We say $(u_1,\ldots,u_t)$ is $(n_1,\ldots,n_t)$-joined 
to $K$ in $G$ by a set of paths $\mathcal{P}$ if the paths satisfy the following properties:
\begin{enumerate}
\item
Every path in $\mathcal{P}$ is a $\{u_1,\ldots,u_t\}$--$V(K)$  path.
\item
$\mathcal{P}$ is a set of internally disjoint paths.
\item
No two paths in $\mathcal{P}$ have the same pair of endpoints.
\item
Exactly $n_i$ paths in $\mathcal{P}$ have $u_i$ as an endpoint, for all $1 \le i \le t$.
\end{enumerate}

\section{Unavoidable configurations}

The main technique  used in this paper is essentially the same as used by Mader
in \cite{M2}. We consider ordered pairs of the form $(G,K)$, where $K$ is an ordered
clique in the graph $G$. We define a reduction operation on such pairs. 

Let $G$ be a graph and $K$ an ordered clique in $G$, such that $V(K) \subset 
V(G)$. Let  $v_1 < v_2 < \cdots < v_k$ be the ordering of the vertices in $K$.
\begin{enumerate}
\item
Suppose there exists a vertex $w \in V(G)\setminus V(K)$ that is adjacent to all 
vertices in $V(K)$. Let $K'$ be the ordered clique in $G$ with  $V(K') = V(K) \cup 
\{w\}$ and the ordering $w < v_1 < \cdots < v_k$ of $V(K')$. We say the
pair $(G,K')$ is obtained from the pair $(G,K)$ by adding the vertex $w$. 
\item
Suppose every vertex in $V(G)\setminus V(K)$ is not adjacent to at least one
vertex in $V(K)$. For every vertex $u \in N_G(v_1) \setminus V(K)$ let $f(u)$ be 
the smallest index such that $v_{f(u)} \not\in N_G(u)$. Let $G'$ be the graph 
obtained from $G$ by deleting the vertex $v_1$ and adding the edge $uv_{f(u)}$, 
for all vertices $u \in N_G(v_1) \setminus V(K)$. Let $K'$ be the ordered clique in
$G'$ with $V(K') = V(K) \setminus \{v_1\}$ and the ordering $v_2 < \cdots < v_k$ 
of $V(K')$. We say the pair $(G',K')$ is obtained from $(G,K)$ by deleting the 
vertex $v_1$.
\end{enumerate}

Note that for any pair $(G,K)$ with $V(K) \subset V(G)$, exactly one of the two 
operations can be applied.  We say a pair $(G',K')$ can be derived from the pair 
$(G,K)$, denoted $(G,K) \rightarrow (G',K')$, if $(G',K')$ can be obtained from 
$(G,K)$ by a sequence of vertex deletion or addition operations.

\begin{lem}
\label{red}
If $(G,K) \rightarrow (G',K')$ then the following properties hold.
\begin{enumerate}
\item
$V(G') \setminus V(K') \subseteq V(G) \setminus V(K)$. 
\item
For every vertex $u \in V(G') \setminus V(K')$, $d_{G'}(u) = d_G(u)$.
\item
For any subset $S$ $\subseteq$ $V(G'-V(K'))$ $\cup$ $E(G'-V(K'))$,\\ $(G-S,K)$ $\rightarrow$ $(G'-S,K')$.
\item 
$|V(G')| + |V(G') \setminus V(K')|$ $<$ $|V(G)| + |V(G) \setminus V(K)|$.
\end{enumerate}
\end{lem}

\pf
The proof follows by induction on the number of reduction operations used
to derive $(G',K')$ from $(G,K)$. It is easy to check that each operation satisfies
the required properties.
\qed

\begin{lem}
\label{del}
Suppose $(G',K')$ is obtained from $(G,K)$ by deleting a vertex. If $(u_1,\ldots,u_t)$ 
is $(n_1,\ldots,n_t)$-joined to $K'$ in $G'$, then it is also $(n_1,\ldots,n_t)$-joined to 
$K$ in $G$.
\end{lem}

\pf Let $K$ be the ordered clique $v_1 < v_2 < \cdots < v_k$ and suppose $K'$ is 
obtained from $K$ by deleting $v_1$, keeping the order of the remaining vertices 
the same. Let $\mathcal{P}$ be the set of paths that $(n_1,\ldots,n_t)$-joins 
$(u_1,\ldots,u_t)$ to $K'$ in $G'$. By the definition of the reduction operation, 
the only edges in $G'$ that are not in $G$ are edges of the form $uv_{f(u)}$ for 
every vertex $u \in N_G(v_1) \setminus V(K)$.  We call such edges 
\emph{bad} edges. Note that there is at most one bad edge incident with any 
vertex $u \in  V(G')\setminus V(K')$, and it must have one endpoint in $V(K')$. Also,
if $uv_{f(u)}$ is a bad edge, by the definition of $f(u)$, $uv_i$ is an edge in $G$, 
for all $1 \le i < f(u)$.

Let $\mathcal{P}_i$ be the set of $n_i$ paths in $\mathcal{P}$ that form a 
$u_i$--$V(K')$ fan.  If none of these paths  contain a bad edge, these form a $u_i$--$V(K)$
fan in $G$.  If any of these paths contains a bad edge, it must be the last edge in the path. 
Let $w_mv_{j_m}$, for $1 \le m \le l$,  be the bad edges contained in the paths in 
$\mathcal{P}_i$, where $1 = j_0  < j_1 < j_2 < \cdots < j_l$,  and $1 \le l \le n_i$.  
Then, replacing the bad edge $w_mv_{j_m}$ by the edge $w_mv_{j_{m-1}}$, 
for $1 \le m \le l$, gives a set of $n_i$ paths that form a $u_i$--$V(K)$ fan in $G$.
These paths have the same set of internal vertices as the paths in  $\mathcal{P}_i$. 
Since the paths in $\mathcal{P}_i$ and $\mathcal{P}_j$ can only have vertices in $V(K')$ 
in common for $i \neq j$, this step can be done independently for each $u_i$. Thus we 
get a set of paths that $(n_1,\ldots,n_t)$-joins $(u_1,\ldots,u_t)$ to $K$ in $G$.
\qed

\begin{lem}
\label{joinone}
Suppose $(G,K')$ is obtained from $(G,K)$ by adding a vertex $w$. Suppose
 $(u_1,\ldots,u_t)$ is $(n_1,\ldots,n_t)$-joined to $K'$ in $G$ by a set of paths 
$\mathcal{P}$. If the sequence $n_1,\ldots,n_t$ does not have a unique maximum, 
and at most one path in $\mathcal{P}$ has $w$ as an endpoint, then 
$(u_1,\ldots,u_t)$ is $(n_1,\ldots,n_t)$-joined to $K$ in $G$.
\end{lem}

\pf
If none of the paths in $\mathcal{P}$ has $w$ as an endpoint, then $\mathcal{P}$ 
is a set of paths that $(n_1,\ldots,n_t)$-joins $(u_1,\ldots,u_t)$ to $K = K'-w$ in
$G$. Suppose $\mathcal{P}$ contains exactly one path terminating in $w$. Without 
loss of generality, we may assume $u_1$ is its other endpoint. Since the sequence 
$n_1,\ldots,n_t$ does not have a unique maximum, there exists an 
$i > 1$ such that $n_i \ge n_1$. Without loss of generality, assume $i = 2$. 
Since $\mathcal{P}$ contains $n_2$ $u_2$--$V(K')$ paths having distinct endpoints 
in $V(K') \setminus \{w\}$, we must have $|V(K')| > n_2 \ge n_1$. 
This implies $|V(K)| \ge n_1$.  Since $\mathcal{P}$ contains $n_1$ $u_1$--$V(K')$ 
paths, one of which terminates in $w$, there exists a vertex $v \in V(K')\setminus 
\{w\}$, such that there is no $u_1$--$v$ path in $\mathcal{P}$. Since $w$ is adjacent 
to all vertices in $K$, adding the edge $wv$ to the $u_1$--$w$ path in $\mathcal{P}$, 
together with all other paths in $\mathcal{P}$, gives a set of paths that 
$(n_1,\ldots,n_t)$-joins $(u_1,\ldots,u_t)$ to $K$ in $G$.
\qed

In view of Lemma~\ref{del}, we will henceforth only need to
consider cases where $(G,K')$ is obtained from $(G,K)$ by adding a vertex $w$.
Suppose $(u_1,\ldots,u_t)$ is $(n_1,\ldots,n_t)$-joined to $K'$ in $G$ by a set of
paths $\mathcal{P}$. In any such case, we will denote by $\mathcal{P}' \subseteq 
\mathcal{P}$ the subset of paths that terminate in $w$, and by $U' \subseteq 
\{u_1,\ldots,u_t\}$ the endpoints of paths in $\mathcal{P}'$ other than $w$. 
We will only consider cases where the sequence $n_1,\ldots,n_t$ does not have a 
unique maximum, and hence by Lemma~\ref{joinone}, we only need to consider 
cases where $|\mathcal{P}'| =|U'| \ge 2$. 

\begin{lem}
\label{jointwo}
Suppose $(G,K')$ is obtained from $(G,K)$ by adding a vertex $w$. Suppose 
$(u_1,\ldots,u_t)$ is $(n_1,\ldots,n_t)$-joined to $K'$ in $G$ by a set of paths 
$\mathcal{P}$ and suppose $|\mathcal{P}'| \ge 2$. Then $(u_1,u_2,\ldots,u_t,w)$ is 
$(n_1',\ldots,n_t',m)$-joined to $K$ in $G-I(\mathcal{P}')$, where $n_i' = n_i-1$ if 
$u_i \in U'$ else $n_i' = n_i$, and $m = \max_{1 \le i \le t}\ n_i'$.
\end{lem}

\pf
Since the paths in $\mathcal{P}$ are internally disjoint, $\mathcal{P} \setminus 
\mathcal{P}'$ is a set of paths that $(n_1',\ldots,n_t')$-joins 
$(u_1,\ldots,u_t)$ to $K$ in $G-I(\mathcal{P}')$. Since any two paths in $\mathcal{P}$ 
can have at most one endpoint in common, $|V(K)| \ge m = \max_{1 \le i \le t}\ n_i'$. 
Since $w$ is adjacent to every vertex in $V(K)$, adding $m$ edges joining $w$ to $V(K)$ 
to the set of paths $\mathcal{P} \setminus \mathcal{P}'$ gives the required set of paths
that $(n_1',\ldots,n_t',m)$-joins $(u_1,\ldots,u_t,w)$ to $K$ in $G-I(\mathcal{P}')$.
\qed

Let $\mathcal{C}$ be a set of graphs such that $A = \{a_1,\ldots,a_t\}$  $\subseteq 
V(H)$, for all graphs $H \in \mathcal{C}$.  Suppose each vertex $a_i \in A$ is assigned 
a positive integer weight $n_i$, for $1 \le i \le t$.  We call such a set of graphs 
$\mathcal{C}$ a \emph{configuration} with \emph{terminal} vertices $(a_1,\ldots,a_t)$
having weights $(n_1,\ldots,n_t)$.

Let $\mathcal{C}$ be a configuration with terminal vertices $(a_1,\ldots,a_t)$ having
weights $(n_1,\ldots,n_t)$. We say that $\mathcal{C}$ is \emph{unavoidable}
if for every graph $G$ and $(G',K')$ such that $(G,\emptyset) \rightarrow (G',K')$,
the following property holds.
\begin{itemize}
\item
If  $(u_1, \ldots, u_t)$  is $(n_1, \ldots, n_t)$-joined to $K'$ in $G'$, then $G$ 
contains a subdivision of some graph $H \in \mathcal{C}$ such that the vertex 
$u_i$ in $G$ corresponds to the vertex $a_i$ in $H$, for $1 \le i \le t$.
\end{itemize}

The basic idea to prove that a configuration $\mathcal{C}$ is unavoidable is to
use induction on the length of the sequence of reductions $(G, \emptyset) =
(G_0,K_0)$, $(G_1,K_1), \ldots, (G_l,K_l)$ such that $(G_{i+1},K_{i+1})$
is obtained from $(G_i,K_i)$ by addition or deletion of vertices. If $\mathcal{C}$
has $t$ terminals $(a_1,\ldots,a_t)$ of weights $(n_1,\ldots,n_t)$, we assume
$(u_1,\ldots,u_t)$ is $(n_1,\ldots,n_t)$-joined to $K_l$ in $G_l$.  In all configurations
that we consider, the sequence of weights does not have a unique maximum.
Lemmas~\ref{del} and \ref{joinone} then imply that if $(G_l,K_l)$ is
obtained from $(G_{l-1},K_{l-1})$ by deleting a vertex or if $|\mathcal{P}'| = 1$,
we can apply induction. If $|\mathcal{P}'| \ge 2$, we use Lemma~\ref{jointwo}
and an appropriate configuration $\mathcal{C}'$ that is either known or assumed to
be unavoidable as part of the induction hypothesis, and apply induction. This may
require that several configurations are proved unavoidable simultaneously.

The following lemma gives a starting point for applying this argument to a graph
of minimum degree at least $d$.

\begin{lem}
\label{edge}
Let $G$ be a graph of minimum degree at least $d \ge 2$.  Then there exists a pair
$(G',K')$ such that $(G,\emptyset) \rightarrow (G',K')$ and $G'-V(K')$ contains an
edge $u_1u_2$ such that $(u_1,u_2)$ is $(d-1,d-1)$-joined to $K'$ in $G'$.
\end{lem}

\pf
Let $(G,\emptyset) = (G_0,K_0)$, $(G_1,K_1), \ldots, (G_l,K_l)$ be
a maximal sequence of pairs such that $(G_i,K_i)$ is obtained from $(G_{i-1},K_{i-1})$
by either deleting or adding a vertex, for $1 \le i \le l$.  Such a sequence exists
since $|V(G_{i+1})|+|V(G_{i+1}) \setminus V(K_{i+1})| < |V(G_i)|+|V(G_i) \setminus 
V(K_i)|$. Then we must have $V(G_l)
= V(K_l)$, otherwise we can add one more pair to the sequence.  Let $i < l$ be the largest
index such that $G_i-V(K_i)$ contains an edge $u_1u_2$.  The choice of $i$ implies
that $(G_{i+1},K_{i+1})$ is obtained from $(G_i,K_i)$ be adding either the vertex $u_1$
or $u_2$ to $K_i$, otherwise $u_1u_2$ is an edge in $G_{i+1}-V(K_{i+1})$.  Without
loss of generality, $V(K_{i+1}) = V(K_i) \cup \{u_1\}$. Then $u_2$ cannot be adjacent
to any vertex other than $u_1$ in $G_i-V(K_i)$.  Since $G$ has minimum degree at
least $d$, $u_2$ has at least $d-1$ neighbors in $V(K_i)$ and thus $|V(K_i)| \ge d-1$. 
Since $u_1$ is adjacent to every vertex in $V(K_i)$, it has at least $d-1$ neighbors in
$V(K_i)$.  Thus $(u_1,u_2)$ is $(d-1,d-1)$-joined to $K_i$ in $G_i$, and $(G_i,K_i)$
is the required pair.  
\qed

Let $H$ be any graph of order $d+1$ and $\mathcal{C}(H)$ the configuration containing 
all possible graphs $H-a_1a_2$, for every edge $a_1a_2 \in E(H)$, with terminal vertices 
$(a_1,a_2)$ having weights $(d-1,d-1)$. If this configuration is unavoidable, 
Lemma~\ref{edge} implies that $H$ is good.

We illustrate the method by restating the proof of Mader's theorem in terms of
unavoidable configurations.  Let $\mathcal{C}(d)$ be the configuration containing the 
single graph $K_{2,d}$, with the two vertices in the part of size 2 being the terminal 
vertices having weight $d$. We claim that for all $d \ge 1$, the configuration 
$\mathcal{C}(d)$ is unavoidable. 

Applying the general strategy,  we may assume $(G_l,K_l)$ is obtained from 
$(G_{l-1},K_{l-1})$ by adding a vertex $w$ and $U' = \{u_1,u_2\}$. Then the union 
of the two paths in  $\mathcal{P}'$ is a $u_1$--$u_2$ path $P$  in $G_{l-1}-V(K_{l-1})$ 
that contains $w$.  If $d = 1$, this gives a subdivision of $K_{2,1}$ in $G$, otherwise  by 
Lemma~\ref{jointwo}, $(u_1,u_2)$ is $(d-1,d-1)$-joined to $K_{l-1}$ in $G_{l-1}- I(P)$. 
By induction, $G-I(P)$ contains  a subdivision of $K_{2,d-1}$ with vertices $u_1,u_2$ 
corresponding to the two terminals in $K_{2,d-1}$. The union of this with the path $P$ 
gives the required subdivision of $K_{2,d}$.  The unavoidability of $\mathcal{C}(d)$ and 
Lemma~\ref{edge} proves Mader's theorem for $d \ge 2$.

Turner's theorem for wheels can be proved in a similar way. In this case, we consider
the configuration $W_d - a_1a_2$, where $a_1a_2$ is a spoke and $a_1$ the center of 
the wheel. Both $a_1,a_2$ have weight $d-1$.  We also need another configuration $W_d -
\{a_1a_2, a_1a_3, a_2a_3\}$, where $\{a_1,a_2,a_3\}$ induce a triangle in $W_d$
with $a_1$ being the center of the wheel. If $d = 3$, $a_1,a_2,a_3$ all have weight 1,
while for $d \ge 4$, $a_1,a_2$ have weight $d-2$ and $a_3$ has weight $d-3$. It
can be argued in a similar way that both these configurations are unavoidable for
all $d \ge 3$.

\section{Maximal 3-degenerate graphs}

A maximal 3-degenerate graph of order $n \ge 3$ is a graph whose vertices 
can be ordered  $v_1,\ldots,v_n$ such that $\{v_1,v_2,v_3\}$ induce a $K_3$ and $v_i$
is adjacent to exactly 3 vertices in $\{v_1,\ldots,v_{i-1}\}$, for $4 \le i \le n$.

\begin{thm}
\label{3-deg}
Every graph of minimum degree at least $d \ge 2$ contains a subdivision of
some maximal 3-degenerate graph of order $d+1$.
\end{thm}

\pf
The proof follows the same general strategy. We define a set of configurations 
and show that they are unavoidable. The theorem then follows by applying
Lemma~\ref{edge}.

Consider the following configurations.
\begin{enumerate}
\item
$\mathcal{C}_1(d)$ for $d \ge 1$, contains all graphs of order $d+2$ with two
terminal vertices $a_1,a_2$, and $d$ other vertices $b_1, \ldots,  b_d$, ordered
so that $b_1$ is adjacent to $a_1$ and $a_2$, and $b_i$ is adjacent to exactly
3 vertices in $\{a_1,a_2,b_1,\ldots,b_{i-1}\}$, for $2 \le i \le d$. The two terminal 
vertices $a_1,a_2$ have weight $d$.
\item
$\mathcal{C}_2(d)$ for $d \ge 1$, contains all graphs of order $d+3$ with 3
terminal vertices $a_1,a_2,a_3$, and $d$ other vertices $b_1,b_2,\ldots,b_d$, 
ordered so that $b_i$ is adjacent to exactly 3 vertices in 
$\{a_1,a_2,a_3,b_1,\ldots,b_{i-1}\}$, for $1 \le i \le d$. The 3 terminal vertices have 
weight $d$ each.
\item
$\mathcal{C}_3(d)$ for $d \ge 1$, contains all graphs of order $d+3$ with 3
terminal vertices $a_1,a_2,a_3$, such that $a_1$ adjacent to $a_2$, and $d$ other vertices 
$b_1,\ldots,b_d$, ordered so that $b_i$ is adjacent to exactly 3 vertices in 
$\{a_1,a_2,a_3,b_1, \ldots,b_{i-1}\}$ for $1 \le i \le d$. The vertices $a_1,a_2$ have 
weight $d+1$, while $a_3$ has weight $d$.
\item
$\mathcal{C}_4(d,t)$ for $d \ge 0$, $t \ge 3$, contains all graphs with $t$ terminal vertices
$a_1,a_2,\ldots,a_t$ such that $a_1$ is adjacent to $a_i$ for $2 \le i \le t$, and $d$ other
vertices $b_1,\ldots,b_d$ such that $b_i$ is adjacent to exactly 3 vertices in
$\{a_1,\ldots,a_t,b_1,\ldots,b_{i-1}\}$ for $1 \le i \le d$. The weight of $a_1$ is $d+t-1$,
and the weight of $a_i$ is $d+i-1$, for $2 \le i \le t$.
\end{enumerate}

We show that the configurations $\mathcal{C}_1$, $\mathcal{C}_2$, $\mathcal{C}_3$, 
$\mathcal{C}_4$ are unavoidable. We consider each of the 4 configurations.

\noindent {\bf Case 1.} Consider the configuration $\mathcal{C}_1$.  If $d=1$, this just 
contains the graph $K_{2,1}$ with 2 terminal vertices of weight 1. This is unavoidable, as 
argued in the proof of Mader's theorem. Suppose $d \ge 2$.   We may assume 
$\mathcal{P}'$ contains exactly 2 paths. Lemma~\ref{jointwo} implies $(u_1,u_2,w)$ is 
$(d-1,d-1,d-1)$-joined to $K_{l-1}$ in $G_{l-1}-I(\mathcal{P}')$. By induction, 
$G-I(\mathcal{P}')$ contains a subdivision of some graph in $\mathcal{C}_2(d-1)$, 
with vertices $u_1,u_2,w$  corresponding to $a_1,a_2,a_3$ , respectively.  Adding the 
paths in $\mathcal{P}'$ to this gives a subdivision of a graph in $\mathcal{C}_1(d)$.

\noindent {\bf Case 2.} Consider the configuration $\mathcal{C}_2$. If $d=1$, this 
contains the graph $K_{3,1}$ with 3 terminal vertices of weight 1. If  $\mathcal{P}'$ 
contains 3 paths, this gives a subdivision of $K_{3,1}$ in $G_{l-1}-V(K_{l-1})$, with 
$u_1,u_2,u_3$ corresponding to $a_1,a_2,a_3$, respectively.
Suppose $|\mathcal{P}'| = 2$, and assume without loss of generality $U' = \{u_1,u_2\}$. 
Lemma~\ref{jointwo} implies $(u_3,w)$ is $(1,1)$-joined to $K_{l-1}$ in $G_{l-1}-
I(\mathcal{P}')$. This implies $G-I(\mathcal{P}')$ contains a $u_3$--$w$ path.
Adding this to the paths in $\mathcal{P}'$ gives a subdivision of $K_{3,1}$ in 
which $u_1,u_2,u_3$ correspond to $a_1,a_2,a_3$, respectively. 

A similar argument holds if $d \ge 2$. If $\mathcal{P}'$ contains 3 paths, by 
Lemma~\ref{jointwo}, $(u_1,u_2,u_3)$ is $(d-1,d-1,d-1)$-joined to $K_{l-1}$ in  
$G_{l-1}-(I(\mathcal{P}') \cup \{w\})$. By induction, $G-(I(\mathcal{P}') \cup \{w\})$ 
contains a subdivision of some graph in $\mathcal{C}_2(d-1)$, with vertices 
$u_1,u_2,u_3$ corresponding to $a_1,a_2,a_3$, respectively.  Adding the
vertex $w$ and the paths in $\mathcal{P}'$ to this, gives a subdivision of a graph in 
$\mathcal{C}_2(d)$. Suppose 
$|\mathcal{P}'| = 2$ and assume without loss of generality, $U' = \{u_1,u_2\}$.
Lemma~\ref{jointwo} implies $(u_3,w,u_1)$ is $(d,d,d-1)$-joined to $K_{l-1}$ in
$G_{l-1}-I(\mathcal{P}')$. Therefore $G-I(\mathcal{P}')$ contains a subdivision of some 
graph in $\mathcal{C}_3(d-1)$, with vertices $u_3,w,u_1$  corresponding to 
$a_1,a_2,a_3$, respectively. Adding the paths in $\mathcal{P}'$ to this gives a subdivision 
of a graph in $\mathcal{C}_2(d)$.

\noindent {\bf Case 3.} Consider the configuration $\mathcal{C}_3$.  Suppose 
$|\mathcal{P}'|=3$. If $d = 1$, then $(u_1,u_2)$ is $(1,1)$-joined to $K_{l-1}$ in
$G_{l-1}- (I(\mathcal{P}') \cup \{w\})$.  This implies $G-(I(\mathcal{P}') \cup \{w\})$ 
contains a $u_1$--$u_2$ path. Adding $w$ and the paths in $\mathcal{P}'$ to this, 
gives the required subdivision of the graph in $\mathcal{C}_3(1)$.  If $d \ge 2$, then 
$(u_1,u_2,u_3)$ is $(d,d,d-1)$-joined to $K_{l-1}$ in $G_{l-1}-(I(\mathcal{P}') \cup 
\{w\})$. By induction, $G-(I(\mathcal{P}') \cup \{w\})$ contains a subdivision of some 
graph in $\mathcal{C}_3(d-1)$, with vertices $u_1,u_2,u_3$  corresponding to 
$a_1,a_2,a_3$, respectively. Adding the vertex $w$ to this
along with the paths in $\mathcal{P}'$, gives the required subdivision of a graph in
$\mathcal{C}_3(d)$.

Suppose $|\mathcal{P}'| = 2$ and $U' = \{u_1,u_2\}$. 
The union of the two paths in $\mathcal{P}'$ is a $u_1$--$u_2$ path $P$ in 
$G_{l-1}-V(K_{l-1})$ that contains $w$. Lemma~\ref{jointwo} implies $(u_1,u_2,u_3)$ is 
$(d,d,d)$-joined to $K_{l-1}$ in $G_{l-1}-I(P)$. By induction, $G-I(P)$ contains a 
subdivision of some graph in $\mathcal{C}_2(d)$, with vertices $u_1,u_2,u_3$ 
corresponding to $a_1,a_2,a_3$, respectively.  Adding the path $P$ to this gives a 
subdivision of some graph in $\mathcal{C}_3(d)$.

Suppose $|\mathcal{P}'| = 2$ and $U' = \{u_2,u_3\}$. The
case when $U' = \{u_1,u_3\}$ is symmetric. Then $(u_1,u_2,w)$ is 
$(d+1,d,d+1)$-joined to $K_{l-1}$ in $G_{l-1}-I(\mathcal{P}')$. By induction, $G-
I(\mathcal{P}')$ contains a subdivision of some graph in $\mathcal{C}_4(d-1,3)$, with 
vertices $u_1,u_2,w$  corresponding to $a_1,a_2,a_3$, respectively.  Adding the paths in 
$\mathcal{P}'$ to it, gives a subdivision of a graph in $\mathcal{C}_3(d)$.

\noindent {\bf Case 4.} Consider the configuration $\mathcal{C}_4$. 

\noindent {\bf Case 4.1} Suppose $d = 0$. The only graph in $\mathcal{C}_4(0,t)$ has 
$t \ge 3$ terminals $a_1,\ldots,a_t$ with edges $a_1a_i$, for $2 \le i \le t$.  The weight 
of $a_1$ is $t-1$ and that of $a_i$ is $i-1$ for $2 \le i \le t$.  In this case, we need to 
show that there exist $t-1$ paths in $G$ that form a $u_1$--$\{u_2,\ldots,u_t\}$ fan. 

Suppose $u_1 \in U'$. Let $i$ be the smallest index greater than 1 such that 
$u_i \in U'$. The union of the $u_1$--$w$ and $u_i$--$w$ paths in $\mathcal{P}'$
is a $u_1$--$u_i$ path $P$ in $G_{l-1}-V(K_{l-1})$. If $t > 3$, then 
$(u_1,u_2,\ldots,u_{i-1}, u_{i+1},\ldots,u_t)$ is $(t-2,1,\ldots,i-2,i-1,\ldots,t-2)$-joined 
to $K_{l-1}$ in $G_{l-1}-(I(\mathcal{P}') \cup  \{w\})$. By induction, 
$G-(I(\mathcal{P}') \cup \{w\})$ contains $t-2$ internally disjoint paths that
form a $u_1$--$\{u_2,\ldots,u_{i-1},u_{i+1},\ldots,u_t\}$ fan.  Adding the path 
$P$ to this gives the required set of $t-1$ paths. If $t= 3$, then $(u_1,u_{5-i})$ is
$(1,1)$-joined to $K_{l-1}$ in $G_{l-1}-(I(\mathcal{P}') \cup \{w\})$. Thus $G-
(I(\mathcal{P}') \cup \{w\})$ contains a $u_1$--$u_{5-i}$ path. Adding the path $P$ 
to this gives the required paths that form a $u_1$--$\{u_2,u_3\}$ fan.

Suppose $u_1 \not\in  U'$. Again, let $i$ be the smallest index such that $u_i \in U'$. 
Then $(u_1,\ldots,u_{i-1},u_{i+1},\ldots,u_t,w)$ is 
$(t-1,1,2,\ldots,t-2,t-1)$-joined to $K_{l-1}$ in $G_{l-1}-I(\mathcal{P}')$. By induction,  
$G-I(\mathcal{P}')$ contains $t-1$ internally disjoint paths
that form a $u_1$--$\{u_1,\ldots,u_{i-1}$, $u_{i+1},\ldots,u_t,w\}$ fan. The union of 
the $u_1$--$w$ path in this set with the  $u_i$--$w$ path in $\mathcal{P}'$ is a 
$u_1$--$u_i$ path in $G$ that is internally disjoint from the other paths in the set. 
Replacing the $u_1$--$w$ path in the fan by this gives $t-1$ paths that form a
$u_1$--$\{u_2,\ldots,u_t\}$ fan.

\noindent {\bf Case 4.2} Suppose $d > 0$. If $|\mathcal{P}'| \ge 3$, then 
$(u_1,\ldots,u_t)$ is $(d+t-2, d,\ldots,d+t-2)$-joined to $K_{l-1}$ in $G_{l-1}-
(I(\mathcal{P}') \cup \{w\})$.
By induction, $G-(I(\mathcal{P}')\cup \{w\})$ contains a subdivision of some graph in
$\mathcal{C}_4(d-1,t)$ with vertex $u_i$ corresponding to $a_i$, for $1 \le i \le t$. Adding
the vertex $w$ and any 3 paths in $\mathcal{P}'$, we get a subdivision of graph in 
$\mathcal{C}_4(d,t)$ that is contained in $G$. 

Suppose $|\mathcal{P}'| = 2$ and $u_1 \in U'$. Let $u_i, i > 1$ be the other vertex in 
$U'$. Then the union of the two paths in $\mathcal{P}'$ is  a $u_1$--$u_i$ path $P$ in 
$G_{l-1}-V(K_{l-1})$.  If $t = 3$, then $(u_1,u_{5-i},u_i)$ is $(d+1,d+1,d)$-joined to 
$K_{l-1}$ in $G_{l-1}-I(P)$. By induction, $G-I(P)$ contains a subdivision of some graph 
in $\mathcal{C}_3(d)$ with vertices $u_1,u_{5-i},u_i$  corresponding to $a_1,a_2,a_3$ , 
respectively. Adding the $u_1$--$u_i$ path $P$ to this gives a subdivision of a graph in 
$\mathcal{C}_4(d,3)$. If $t > 3$, then $(u_1,\ldots,u_{i-1},u_{i+1},\ldots,u_t)$ is 
$(d+t-2,d+1,\ldots,d+i-2,d+i-1,\ldots,d+t-2)$-joined to $K_{l-1}$ in $G_{l-1}-I(P)$. 
By induction, $G-I(P)$ contains a subdivision of some graph in $\mathcal{C}_4(d,t-1)$ 
with vertices $u_1,\ldots,u_{i-1},u_{i+1},\ldots,u_t$  corresponding to 
$a_1,\ldots,a_{t-1}$, respectively.
Adding the path $P$ to this gives a subdivision of a graph in $\mathcal{C}_4(d,t)$.

Finally, suppose $|\mathcal{P}'| = 2$ and $u_1 \not\in U'$. Then 
$(u_1, u_2,\ldots,u_t,w)$ is $(d+t-1, d, \ldots, d+t-2,d+t-1)$-joined to
$K_{l-1}$ in $G_{l-1}-I(\mathcal{P}')$. By induction, $G-I(\mathcal{P}')$ contains a 
subdivision of some graph in $\mathcal{C}_4(d-1,t+1)$, with vertices 
$u_1,\ldots,u_t,w$ corresponding to $a_1,\ldots,a_{t+1}$, respectively. Adding the 
paths in $\mathcal{P}'$ to this gives a subdivision of a graph in $\mathcal{C}_4(d,t)$.

This completes all cases and we conclude that all the 4 configurations are unavoidable.
The theorem then follows from Lemma~\ref{edge} and the fact that the configuration
$\mathcal{C}_1(d-1)$ is unavoidable. Note that  for any graph in $\mathcal{C}_1(d-1)$, 
adding an edge between the two terminal vertices $a_1,a_2$ gives a maximal 3-degenerate 
graph of order $d+1$.
\qed

\section{Planar Maximal 3-degenerate Graphs}

Since there exist non-planar 3-degenerate graphs of order 6, not every maximal
3-degenerate graph is good.  However, we do not know of any planar graph
that is not good. This suggests the following problem.
 
\begin{conj}
\label{planar}
 Every planar maximal 3-degenerate graph is good.
\end{conj}

A specific family of planar maximal 3-degenerate graphs is $P_n^3$ with vertices
$v_1,\ldots,v_n$ and $v_i$ adjacent to $v_j$ iff $1 \le |j-i| \le 3$.
We verify  Conjecture~\ref{planar} for two graphs $P_6^3$ and $P_7^3$. Note
that  $P_4^3$ is $K_4$, $P_5^3$ is $K_5^-$ and $P_6^3$ is the only planar 
maximal 3-degenerate graph of order 6.

\begin{thm}
\label{P6}
Every graph of minimum degree at least 5 contains a subdivision of $P_6^3$.
\end{thm}

\pf
The proof is again based on the same technique, using more restricted configurations
than those used in Theorem~\ref{3-deg}. Consider the following set of
configurations.
\begin{enumerate}
\item
$\mathcal{C}_5$ contains a subset of the graphs in the configuration $\mathcal{C}_1(4)$.
The graphs have 6 vertices $\{a_1,a_2,b_1,b_2,b_3,b_4\}$, where $a_1,a_2$ are
terminal vertices of weight 4. The edge sets of the 3 graphs are
\begin{enumerate}
\item
$\{a_1b_1, a_2b_1, a_1b_2,a_2b_2, b_1b_2, a_1b_3, a_2b_3, b_2b_3, a_1b_4,
b_2b_4, b_3b_4\}$.
\item
$\{a_1b_1, a_2b_1, a_1b_2,a_2b_2, b_1b_2, a_1b_3, b_1b_3, b_2b_3, b_1b_4,
b_2b_4, b_3b_4\}$.
\item
$\{a_1b_1, a_2b_1, a_1b_2,a_2b_2, b_1b_2, a_2b_3, b_1b_3, b_2b_3, a_2b_4,
b_2b_4, b_3b_4\}$.
\end{enumerate}
\item
$\mathcal{C}_6$ contains a subset of the graphs in the configuration $\mathcal{C}_2(3)$.
The graphs have 6 vertices $\{a_1,a_2,a_3,b_1,b_2,b_3\}$, where $a_1,a_2,a_3$
are terminal vertices of weight 3. The edge sets of the 3 graphs are
\begin{enumerate}
\item
$\{a_1b_1,a_2b_1,a_3b_1, a_1b_2, a_2b_2, b_1b_2, a_1b_3,b_1b_3,b_2b_3\}$.
\item
$\{a_1b_1,a_2b_1,a_3b_1, a_1b_2, a_3b_2, b_1b_2, a_3b_3, b_1b_3, b_2b_3\}$.
\item
$\{a_1b_1, a_2b_1,a_3b_1, a_2b_2,a_3b_2,b_1b_2, a_2b_3,b_1b_3,b_2b_3\}$.
\end{enumerate}
\item
$\mathcal{C}_7$ contains only one graph from the configuration $\mathcal{C}_2(2)$.
This graph has 5 vertices $\{a_1,a_2,a_3,b_1,b_2\}$, where $a_1,a_2,a_3$ are
terminal vertices of weight 2.  The edges in the graph are $\{a_1b_1$, $a_2b_1$,
$a_3b_1$, $a_1b_2$, $a_2b_2$, $b_1b_2\}$.
\item
$\mathcal{C}_8$ contains only one graph from the configuration $\mathcal{C}_3(2)$.
This graph has 5 vertices $\{a_1,a_2,a_3,b_1,b_2\}$, where $a_1,a_2,a_3$ are
terminal vertices, $a_1,a_2$ have weight 3 and $a_3$ has weight 2. The edge set
of the graph is $\{a_1a_2, a_1b_1,a_2b_1,a_3b_1,a_1b_2,a_2b_2,b_1b_2\}$.
\item
$\mathcal{C}_9$ contains two graphs with 5 vertices $\{a_1,a_2,a_3,a_4,b_1\}$,
where $a_1$, $a_2$, $a_3$, $a_4$ are terminals, $a_1,a_2$ have weight 1 and 
$a_3,a_4$ have weight 2. The edge sets of the two graphs are
\begin{enumerate}
\item
$\{a_1b_1,a_2b_1,a_3b_1,a_4b_1\}$.
\item
$\{a_1b_1, a_2b_1, a_3a_4, a_3b_1\}$.
\end{enumerate}
\end{enumerate}

We show that $\mathcal{C}_5$, $\mathcal{C}_6$, $\mathcal{C}_7$,
$\mathcal{C}_8$ and $\mathcal{C}_9$ are unavoidable.

\noindent {\bf Case 1.} Consider the configuration $\mathcal{C}_5$. The only case
to be considered here is if $|\mathcal{P}'| = 2$. Lemma~\ref{jointwo} implies
that $(u_1,u_2,w)$ is $(3,3,3)$-joined to $K_{l-1}$ in $G_{l-1}-I(\mathcal{P}')$.
By induction, $G-I(\mathcal{P}')$ contains a subdivision of one of the graphs (a), (b)
or (c) in $\mathcal{C}_6$, with vertices $u_1,u_2,w$  corresponding
to $a_1,a_2,a_3$, respectively. Adding the paths in $\mathcal{P}'$ to this gives
a subdivision of the corresponding graph (a), (b) or (c) in $\mathcal{C}_5$ with 
$u_1,u_2$ corresponding to $a_1,a_2$ and $w$ corresponding to $b_1$.

\noindent {\bf Case 2.} Consider the configuration $\mathcal{C}_6$. If  
$|\mathcal{P}'| = 3$,
then $(u_1,w,u_2)$ is $(2,2,2)$-joined to $K_{l-1}$ in $G_{l-1}-I(\mathcal{P}')$.
By induction, $G-I(\mathcal{P}')$ contains a subdivision of the graph in $\mathcal{C}_7$,
with vertices $u_1,w,u_2$ corresponding to $a_1,a_2,a_3$, respectively.  Adding the
paths in $\mathcal{P}'$ to this gives a subdivision of the graph (a) in $\mathcal{C}_6$.

Suppose $U' = \{u_1,u_2\}$. The other cases can be argued 
symmetrically.  Then $(u_3,w,u_1)$ is $(3,3,2)$-joined to $K_{l-1}$ in 
$G_{l-1}-I(\mathcal{P}')$. By induction,
$G-I(\mathcal{P}')$ contains a subdivision of the graph in $\mathcal{C}_8$, with 
vertices $u_3,w,u_1$ corresponding to $a_1,a_2,a_3$, respectively. Adding the paths in 
$\mathcal{P}'$ to this gives a subdivision of the graph (b) in $\mathcal{C}_6$. 

\noindent {\bf Case 3.} Consider the configuration $\mathcal{C}_7$.  If 
$|\mathcal{P}'| = 3$, then $(u_1,u_2,w)$ is $(1,1,1)$-joined to $K_{l-1}$ in
$G_{l-1}-I(\mathcal{P}')$.  Since the configuration $\mathcal{C}_2(1)$ is
unavoidable, $G-I(\mathcal{P}')$ contains a subdivision of $K_{3,1}$ in
which the vertices $u_1,u_2,w$ correspond to the vertices in the part of size 3.
Adding the paths in $\mathcal{P}'$ to this, gives a subdivision of the graph
in  $\mathcal{C}_7$.

Suppose $U' = \{u_1,u_2\}$. Then $(u_1,u_2,w,u_3)$ is $(1,1,2,2)$-joined to 
$K_{l-1}$ in $G_{l-1}-I(\mathcal{P}')$ and $G-I(\mathcal{P}')$ contains a subdivision
of one of the two graphs in $\mathcal{C}_9$, with vertices  $u_1,u_2,w,u_3$   
corresponding to $a_1,a_2,a_3,a_4$, respectively. In either case, adding the paths in 
$\mathcal{P}'$ to this gives a subdivision of the graph in $\mathcal{C}_7$. 

Suppose $U' = \{u_2,u_3\}$.  Then $(u_1,w,u_2)$ is $(2,2,1)$-joined to $K_{l-1}$ in
$G_{l-1}-I(\mathcal{P}')$ and $G-I(\mathcal{P}')$ contains a subdivision
of the graph in $\mathcal{C}_3(1)$, with vertices  $u_1,w,u_2$    
corresponding to  $a_1,a_2,a_3$, respectively. Adding the paths in $\mathcal{P}'$
to this gives a subdivision of the graph in $\mathcal{C}_7$. The case when 
$U' = \{u_1,u_3\}$ can be argued symmetrically.

\noindent {\bf Case 4.} Consider the configuration $\mathcal{C}_8$. If $|\mathcal{P}'|
= 3$ then $(u_1,u_2,w)$ is $(2,2,1)$-joined to $K_{l-1}$ in $G_{l-1}-I(\mathcal{P}')$.
By induction, $G-I(\mathcal{P}')$ contains a subdivision of the graph in 
$\mathcal{C}_3(1)$, with vertices  $u_1,u_2,w$   corresponding to $a_1,a_2,a_3$,
respectively.  Together with the paths in $\mathcal{P}'$, this gives a subdivision of
the graph in  $\mathcal{C}_8$.

Suppose $U' = \{u_1,u_2\}$.  The union of the two paths in $\mathcal{P}'$ is a 
$u_1$--$u_2$ path $P$  in $G_{l-1}-V(K_{l-1})$ that contains $w$. Since 
$(u_1,u_2,u_3)$ is
$(2,2,2)$-joined to $K_{l-1}$ in $G_{l-1}-I(P)$, $G-I(P)$ contains a subdivision of
the graph in $\mathcal{C}_7$, with vertices   $u_1,u_2,u_3$  corresponding to
$a_1,a_2,a_3$, respectively. Adding the path $P$ to this gives a subdivision of
the graph in $\mathcal{C}_8$. 

Suppose $U' = \{u_2,u_3\}$.  Then $(u_1,u_2,w)$ is $(3,2,3)$-joined to $K_{l-1}$ in 
$G_{l-1}-I(\mathcal{P}')$.
By induction, $G-I(\mathcal{P}')$ contains a subdivision of the graph in
$\mathcal{C}_4(1,3)$, with vertices  $u_1,u_2,w$  corresponding to  $a_1,a_2,a_3$,
respectively. Adding the paths in $\mathcal{P}'$ to this gives a subdivision of the
graph in $\mathcal{C}_8$. The case when $U' = \{u_1,u_3\}$ is similar.

\noindent {\bf Case 5.} Consider the configuration $C_9$. If $|\mathcal{P}'|=4$, then
$G_{l-1}-V(K_{l-1})$ contains a subdivision of the graph (a) in
$\mathcal{C}_9$, with vertices   $u_1,u_2,u_3,u_4,w$   corresponding to
$a_1,a_2,a_3,a_4,b_1$, respectively.

If $|\mathcal{P}'| = 3$, then $(u_i,w)$ is $(1,1)$-joined to $K_{l-1}$ in
$G_{l-1}-I(\mathcal{P}')$, where $u_i \not\in U'$. 
Then $G-I(\mathcal{P}')$ contains a $u_i$--$w$ path, which together with the
paths in $\mathcal{P}'$  gives a subdivision of the graph (a) in $\mathcal{C}_9$.

Suppose $|\mathcal{P}'|=2$ and $U' \neq \{u_3,u_4\}$. Then $(w,u_i,u_j)$ is 
$(2,1,2)$-joined to $K_{l-1}$ in $G_{l-1}-I(\mathcal{P}')$, where $u_i, u_j \not\in U'$ 
and $1 \le i < j\le 4$.  By induction, $G-I(\mathcal{P}')$ contains a subdivision of the 
graph in $\mathcal{C}_4(0,3)$ with vertices  $w,u_i,u_j$  corresponding to $a_1,a_2,a_3$, 
respectively. Adding the paths in $\mathcal{P}'$ gives a subdivision of the graph (a) in 
$\mathcal{C}_9$.

The only other possibility is that $U' = \{u_3,u_4\}$. The union of the 2 paths in 
$\mathcal{P}'$ is a $u_3$--$u_4$ path $P$ in $G_{l-1}-V(K_{l-1})$. Since 
$(u_1,u_2,u_3)$ is $(1,1,1)$-joined to $K_{l-1}$ in $G_{l-1}-I(P)$, $G-I(P)$ contains a 
subdivision of the graph in $\mathcal{C}_2(1)$, with vertices  $u_1,u_2,u_3$    corresponding to $a_1,a_2,a_3$, respectively. Adding the path $P$ to this gives a 
subdivision of the graph (b) in $\mathcal{C}_9$.

This completes all cases and shows that the configurations $\mathcal{C}_5$,
$\mathcal{C}_6$, $\mathcal{C}_7$, $\mathcal{C}_8$ and $\mathcal{C}_9$ are
unavoidable. Theorem~\ref{P6} then follows from Lemma~\ref{edge}, since
adding the edge $a_1a_2$ to any graph in $\mathcal{C}_5$ gives the graph
$P_6^3$.
\qed

We next consider planar maximal 3-degenerate graphs of order 7. There are 3 different
such graphs, but we consider only the graph $P_7^3$. While it is possible to use
the same technique, the number of configurations required appears to be large. We
can reduce the number of configurations required by starting with an initial graph
other than an edge.

Let $\mathcal{C}$ be a configuration with terminal vertices $a_1,\ldots,a_t$ of
weights $n_1,\ldots$, $n_t$, respectively. We say a pair $(G,K)$ contains the configuration 
$\mathcal{C}$ if $G-V(K)$  contains a subdivision $\mathcal{T}(H)$ of some graph 
$H \in \mathcal{C}$, such that vertices $u_1,\ldots,u_t$ correspond to $a_1,\ldots,a_t$, 
respectively, and  $(u_1,\ldots,u_t)$ is $(n_1,\ldots,n_t)$-joined to $K$ in 
$G-V(\mathcal{T}(H)) \setminus \{u_1,\ldots,u_t\}$. 

Consider the following set of configurations.
\begin{enumerate}
\item
$\mathcal{C}_{10}(d)$ for $d \ge 1$ is the configuration containing only the graph
$K_2$ with 2 terminal vertices of weight $d$.
\item
$\mathcal{C}_{11}(d)$ for $d \ge 1$ is the configuration containing only the graph
$K_3$ with 3 terminal vertices of weight $d$.
\item
$\mathcal{C}_{12}(d)$ for $d \ge 1$ is the configuration containing only the
graph $K_4^-$, obtained by deleting an edge from $K_4$. There are 3 terminal 
vertices $a_1,a_2,a_3$ with $a_1,a_3$ of weight $d+1$ and $a_2$ of weight $d$. 
The missing edge is $a_1a_3$.
\item
$\mathcal{C}_{13}(d)$ for $d \ge 1$ is the configuration containing only the graph
$K_4$ with 3 terminal vertices $a_1,a_2,a_3$ of weight $d$.
\end{enumerate}

\begin{lem}
\label{K4}
Let $G$ be a graph of minimum degree at least $d \ge 4$. Then there exists a pair $(G',K')$
such that $(G,\emptyset) \rightarrow (G',K')$ and  $(G',K')$ contains the
configuration $\mathcal{C}_{13}(d-3)$.
\end{lem}

\pf
Let $(G,\emptyset) = (G_0,K_0)$, $(G_1,K_1), \ldots, (G_l,K_l)$ be a maximal
sequence of pairs such that  $(G_{i+1},K_{i+1})$ is obtained from $(G_i,K_i)$ by
adding or deleting a vertex, for $0 \le i < l$. Let $i$ be the smallest index such that 
$(G_i,K_i)$ contains the configuration $\mathcal{C}_{10}(d-1)$.  Lemma~\ref{edge} 
implies there exists such an index $i$. Since $d \ge 4$, we have $i > 0$.  Then $(G_i,K_i)$ 
must be obtained from $(G_{i-1},K_{i-1})$ by adding a vertex $w$, and 
$U' = \{u_1,u_2\}$. 
This implies $(G_{i-1},K_{i-1})$ contains the configuration $\mathcal{C}_{11}(d-2)$ 
with vertices $u_1,u_2,w$ corresponding to $a_1,a_2,a_3$, respectively.

Let $j$ be the smallest index such that $(G_j,K_j)$ contains $\mathcal{C}_{11}(d-2)$. 
Since $d \ge 4$, we have $j > 0$. Again, $(G_j,K_j)$ must be obtained from 
$(G_{j-1},K_{j-1})$ by adding a vertex $w$. If $|\mathcal{P}'|= 3$, then  
$(G_{j-1},K_{j-1})$ contains $\mathcal{C}_{13}(d-3)$, with vertices $u_1,u_2,u_3$
corresponding to the vertices $a_1,a_2,a_3$, respectively. 
Then $(G_{j-1},K_{j-1})$ is the required pair.

Suppose $|\mathcal{P}'| = 2$, and without loss of generality, $U' = \{u_2,u_3\}$.
Then $(G_{j-1},K_{j-1})$ contains $\mathcal{C}_{12}(d-3)$ with vertices 
$u_1,u_2,w$ corresponding to $a_1,a_2,a_3$, respectively. 

Let $m$ be the smallest index such that $(G_m,K_m)$
contains $\mathcal{C}_{12}(d-3)$. Since $d \ge 4$, we have $m > 0$. Again,
$(G_m,K_m)$ must be obtained from $(G_{m-1},K_{m-1})$ by adding a vertex $w$.
If $|\mathcal{P}'| = 3$, then $(G_{m-1},K_{m-1})$ contains the configuration
$\mathcal{C}_{13}(d-3)$ with vertices $u_1,u_3,w$ corresponding to
$a_1,a_2,a_3$, respectively. If $U' = \{u_1,u_3\}$ then $(G_{m-1},K_{m-1})$ contains 
$\mathcal{C}_{13}(d-3)$ with vertices $u_1,u_2,u_3$ corresponding to $a_1,a_2,a_3$, 
respectively. If $U' = \{u_2,u_3\}$ then $(G_{m-1},K_{m-1})$ contains 
$\mathcal{C}_{12}(d-3)$, with vertices $u_1,u_3,w$ corresponding to
$a_1,a_2,a_3$, respectively. This contradicts the choice of $m$. Similarly,
if $U' = \{u_1,u_2\}$, then $(G_{m-1},K_{m-1})$ contains $\mathcal{C}_{12}(d-3)$, 
with vertices $u_3,u_1,w$ corresponding to $a_1,a_2,a_3$, respectively. Again, 
this contradicts the choice of $m$. Therefore $(G_{m-1},K_{m-1})$ must contain 
$\mathcal{C}_{13}(d-3)$.
\qed

\begin{thm}
\label{P73}
Every graph of minimum degree at least 6 contains a subdivision of $P_7^3$.
\end{thm}

\pf
Let $G$ be a graph of minimum degree at least 6. Lemma~\ref{K4} implies there exists
a pair $(G',K')$ such that $(G,\emptyset) \rightarrow (G',K')$ and $(G',K')$ contains
the configuration $\mathcal{C}_{13}(3)$. Thus $G'-V(K')$ contains a subdivision
$H$ of $K_4$, with vertices $u_1,u_2,u_3$ corresponding to the vertices
$a_1,a_2,a_3$, respectively, such that $(u_1,u_2,u_3)$ is $(3,3,3)$-joined to 
$K'$ in $G'-(V(H)\setminus \{u_1,u_2,u_3\})$. Since $\mathcal{C}_6$ is unavoidable, 
$G-(V(H) \setminus \{u_1,u_2,u_3\})$ contains a subdivision of one of the graphs (a), 
(b),  or (c) in $\mathcal{C}_6$, with vertices $u_1,u_2,u_3$ corresponding to 
$a_1,a_2,a_3$, respectively.  In all cases, the union of this graph with $H$ gives a 
subdivision of $P_7^3$ in $G$.
\qed

\section{Remarks}

We have verified Conjecture~\ref{planar} for the other two planar maximal 
3-degenerate graphs of order 7. Although the method is the same, the number of 
configurations required is larger, and we omit the details. 
A planar maximal 3-degenerate graph is also a
maximal planar graph. An interesting question is whether all maximal planar graphs
are good? The smallest case to consider is the octahedron, obtained by deleting a
perfect matching from $K_6$. While we do not know a graph of minimum degree 5
that does not contain a subdivision of this, the technique used in this paper cannot
be applied since the required configuration is avoidable. It would be interesting
to see if there is any characterization of unavoidable configurations. Perhaps the
first question to answer would be to find the maximum number of edges in
a good graph of order $d+1$. For $2 \le d \le 5$, this is exactly $3d-3$.  Does this
hold in general? An even simpler question would be to find the largest number
$m$ such that every graph of minimum degree $d$ contains a subdivision of
some graph of order $d+1$ and size $m$. Theorem~\ref{3-deg} shows that
$m \ge 3d-3$ and the bound is tight for $2 \le d \le 5$. Does this hold for all d?
Finally, it would be interesting to consider non-separating versions of these results.
Kriesell~\cite{K} generalized Dirac's theorem to show that every connected graph
$G$ with minimum degree at least 4 contains a subdivision $H$ of $K_4$ such that
$G-V(H)$ is connected. Can the results in this paper be extended in a similar way,
by increasing the minimum degree bound by one?

\end{document}